\documentclass[a4paper,12pt]{article}
\usepackage{amsmath,amscd,amssymb}
\begin{document}

\centerline{\textbf{
\textsc{VARI\'ET\'ES ALG\'EBRIQUES DONT LE FIBR\'E
TANGENT}}}
\centerline{\textbf{
\textsc{EST TOTALEMENT D\'ECOMPOS\'E}}}
$\ $
\vspace{1cm}\\
\centerline{St\'ephane \textsc{Druel}}
\begin{center}
DMI-\'Ecole Normale Sup\'erieure\\
45 rue d'Ulm\\
75005 PARIS\\
e-mail: \texttt{druel@clipper.ens.fr}
\end{center}
$\ $
\vspace{1cm}\\
\centerline{\textsc{Introduction}}\\
\newline
\indent Soit $X$ une vari\'et\'e compacte k\"ahl\'erienne dont
le rev\^etement universel $\tilde{X}$ est isomorphe au produit 
$\prod_{i\in I}U_{i}$ de vari\'et\'es complexes lisses et sur lequel le 
groupe $\pi_{1}(X)$ agit diagonalement. La d\'ecomposition 
$T_{\tilde{X}}=\oplus_{i\in I}p_{i}^{*}T_{U_{i}}$ induit alors une 
d\'ecomposition de $T_{X}$ en somme directe de sous-fibr\'es int\'egrables. 
Ce travail contribue \`a l'\'etude de l'assertion r\'eciproque
et compl\`ete les r\'esultats d\'ej\`a obtenus par Beauville ([B1]) :\\
\newline
\textbf{Th\'eor\`eme 1}.$-$\textit{Soient X une vari\'et\'e projective
lisse de dimension $n\ge 1$ dont le fibr\'e tangent est totalement 
d\'ecompos\'e et $T_{X}=M_{1}\oplus\cdots\oplus M_{n}$ ladite 
d\'ecomposition. On suppose que les fibr\'es 
$\oplus_{i\in I}M_{i}$ sont int\'egrables, pour tout ensemble d'indices
$I\subset\{1,\ldots,n\}$. Le rev\^etement universel 
$\tilde{X}$ de $X$ est alors produit de surfaces de Riemann et la
d\'ecomposition de $T_{X}$ est induite 
par la d\'ecomposition canonique de $T_{\tilde{X}}$.}\\
\newline
\indent Rappelons qu'une vari\'et\'e projective lisse est dite
\textit{minimale} si $K_{X}$ est num\'eriquement effectif.\\
\newline
\textbf{Th\'eor\`eme 2}.$-$\textit{Soit X une vari\'et\'e projective
lisse minimale de dimension $n\ge 1$ dont le fibr\'e tangent est totalement 
d\'ecompos\'e. Le rev\^etement universel 
$\tilde{X}$ de $X$ est alors produit de surfaces de Riemann et la
d\'ecomposition de $T_{X}$ est induite 
par la d\'ecomposition canonique de $T_{\tilde{X}}$.}\\
\newline
\indent Nous d\'emontrons ces r\'esultats en exhibant une fibration
lisse de $X$ munie d'une connexion int\'egrable compatible \`a la
d\'ecomposition de $T_{X}$.\\
\newline
\textbf{Remerciements}.$-$Je tiens \`a exprimer toute ma gratitude \`a
Arnaud Beauville pour m'avoir soumis ce probl\`eme et pour l'aide qu'il m'a
apport\'ee.\\
\newline
\centerline{\textsc{D\'emonsration des Th\'eor\`emes}}\\
\newline
\textbf{Lemme 1 }([B1] lemme 3.1).$-$\textit{Soient $X$ une vari\'et\'e 
lisse et E un facteur
direct de $T_{X}$. La classe d'Atiyah 
$at(E)\in H^{1}(X,\Omega_{X}^{1}\otimes\mathcal{E}nd(E))$
provient de $H^{1}(X,E^{*}\otimes\mathcal{E}nd(E))$. En particulier, 
tout \'el\'ement de $H^{r}(X,\Omega_{X}^{r})$ donn\'e par un 
polyn\^ome en les classes de Chern de $E$ est nul, pour $r$ strictement
plus grand que le rang de $E$.}\\
\newline
\textbf{Corollaire 1}.$-$\textit{Soient $X$ une vari\'et\'e projective 
lisse de dimension $n\ge 1$, $E$ un facteur direct de $T_{X}$ de rang 1
et $C\subset X$ une courbe rationnelle irr\'eductible. 
Alors $\text{deg}(E_{|C})=0$ ou $\text{deg}(E_{|C})\ge 2$.}\\
\newline
\textit{D\'emonstration}.$-$Soit $\mathbb{P}^{1}
\overset{v}{\longrightarrow} C$ 
la normalisation de $C$. Supposons $\text{deg}(E_{|C})\le 1$ ;
le groupe de cohomologie 
$H^{1}(\mathbb{P}^{1},E^{*}_{|\mathbb{P}^{1}})$ est donc nul.
Consid\'erons le diagramme commutatif :
\begin{equation*}
\begin{CD}
H^{1}(X,E^{*}) @))) H^{1}(X,\Omega^{1}_{X}) \\
@VVV  @VVV \\
H^{1}(\mathbb{P}^{1},E_{|\mathbb{P}^{1}}^{*}) 
@))) H^{1}(\mathbb{P}^{1},\Omega^{1}_{\mathbb{P}^{1}})
\end{CD}
\end{equation*}
L'\'el\'ement $c_{1}(E)\in H^{1}(X,\Omega_{X}^{1})$ provient de 
$H^{1}(X,E^{*})$ (lemme 1) ; son image 
$c_{1}(E_{|\mathbb{P}^{1}})=\text{deg}(E_{|C})$ 
dans 
$H^{1}(\mathbb{P}^{1},\Omega^{1}_{\mathbb{P}^{1}})$ est 
donc nulle, ce qui prouve le corollaire.\\
\newline
\textbf{Corollaire 2}.$-$\textit{Soient $X$ une vari\'et\'e projective 
lisse de dimension $n\ge 1$ dont le fibr\'e tangent est totalement 
d\'ecompos\'e et $\mathbb{P}^{1}\overset{v}{\longrightarrow} X$
un morphisme non constant. Alors $X$ n'est pas minimale.}\\
\newline
\textit{D\'emonstration}.$-$Ecrivons 
$v^{*}T_{X}=\mathcal{O}_{\mathbb{P}^{1}}(a_{1})
\oplus\cdots\oplus\mathcal{O}_{\mathbb{P}^{1}}(a_{n})$ avec 
$a_{1}\ge\cdots\ge a_{n}\ge 0$ (corollaire 1). L'application tangente 
$T_{\mathbb{P}^{1}}\overset{dv}{\longrightarrow} v^{*}T_{X}$ \'etant 
g\'en\'eriquement injective, on a $a_{1}\ge 2$. On en d\'eduit 
$\text{deg}(v^{*}K_{X})\le -2$ et $X$ n'est donc pas minimale.\\
\newline
\indent Soit $X$ une vari\'et\'e projective lisse complexe. 
Le produit d'intersection entre
1-cycles et diviseurs met en dualit\'e les deux espaces vectoriels
r\'eels $N_{1}(X)=(\{\text{1-cycles}\}/\equiv)\otimes\mathbb{R}$ et 
$N^{1}(X)=(\{\text{diviseurs}\}/\equiv)\otimes\mathbb{R}\,$,
 o\`u $\equiv$ d\'esigne l'\'equivalence num\'erique. 
Soit $NE(X)\subset N_{1}(X)$ le c\^one
engendr\'e par les classes des 1-cycles effectifs. Une \textit{raie
extr\'emale} est une demi-droite $R$ dans $\overline{NE}(X)$, 
adh\'erence de $NE(X)$ dans $N_{1}(X)$, v\'erifiant $K_{X}.R^{*}<0$  et 
telle que pour tout $Z_{1},Z_{2}\in\overline{NE}(X)$, si
$Z_{1}+Z_{2}\in R$ alors $Z_{1},Z_{2}\in R$ ([M2]).
Une \textit{courbe rationnelle extr\'emale} est une 
courbe rationnelle irr\'eductible $C_{0}$ telle que
$\mathbb{R}^{+}[C_{0}]$ soit une raie extr\'emale et telle que
$-K_{X}.C_{0}\le\text{dim}X+1$. 
\textit{Si $X$ est une vari\'et\'e projective 
lisse non minimale alors $X$ contient une courbe rationnelle extr\'emale}
([M2] thm. 1.5). La \textit{longueur} de la raie 
extr\'emale $R$ est ([W]) :
$$\ell(R)=\text{inf}\{-K_{X}.C| C \text{ \'etant une courbe
rationnelle et }C\in R\}.$$
\indent L'\'etude des courbes rationnelles extr\'emales sur les 
vari\'et\'es dont le fibr\'e tangent est totalement 
d\'ecompos\'e fait l'objet du :\\
\newline
\textbf{Lemme 2}.$-$\textit{Soit $X$ une vari\'et\'e projective lisse
non minimale de dimension $n\ge1$ dont le fibr\'e tangent est 
totalement d\'ecompos\'e. Il existe alors un rev\^etement \'etale fini
$Z\longrightarrow X$ tel que $Z$ soit un fibr\'e en droites 
projectives pour la topologie \'etale.}\\
\newline
\textit{D\'emonstration}.$-$ Soit $R=\mathbb{R}^{+}[C_{0}]$ une raie 
extr\'emale engendr\'ee par une courbe rationnelle $C_{0}$ telle 
que $\ell(R)=-K_{X}.C_{0}$. Notons 
$\mathbb{P}^{1}\overset{v_{0}}{\longrightarrow}C_{0}$ 
la normalisation de $C_{0}$.
Les hypoth\`eses fa\^{\i}tes entra\^{\i}nent
la lissit\'e du sch\'ema $\text{Hom}(\mathbb{P}^{1},X)$
(corollaire 1). Soit 
$V\subset\text{Hom}(\mathbb{P}^{1},X)$ la composante connexe
(de dimension $\ell(R)+n$) contenant le point $v_{0}$.
Le groupe $G=PGL_{2}(\mathbb{C})$
agit de mani\`ere naturelle sur $V$ par la formule 
$g.v=vg^{-1}$. Soient $\text{Chow}(X)$ la vari\'et\'e projective
param\'etrant les 1-cycles effectifs et
$V\overset{\alpha}{\longrightarrow}\text{Chow}(X)$ le morphisme
naturel $G$-\'equivariant
qui \`a $\mathbb{P}^{1}\overset{v}{\longrightarrow} X$ associe 
le 1-cycle $v(\mathbb{P}^{1})$ ($v$ est birationnel au dessus de 
$v(\mathbb{P}^{1})$). Enfin, soit $Y$ la normalisation de 
$\overline{\alpha(V)}$ dans le corps $k(V)^{G}$. 
Alors $Y$ est le quotient 
g\'eom\'etrique de $V$ par $G$ 
et l'action de $G$ sur $V$ est libre ([M1] lemme 9 et [W] Appendice A4) ; 
$Y$ est une vari\'et\'e projective et lisse de dimension
$n+\ell(R)-3$. Soit 
$\mathbb{P}^{1}\times V\overset{F}{\longrightarrow} X\times Y$
le morphisme naturel et soit 
$Z=\text{Spec}((F_{*}\mathcal{O}_{\mathbb{P}^{1}\times V})^{G})$.
Alors $Z$ est le quotient g\'eom\'etrique de $\mathbb{P}^{1}\times V$
par $G$, l'action de $G$ \'etant donn\'ee par la formule
$g.(z,v)=(g(z),vg^{-1})$ ; $Z$ est une vari\'et\'e
projective et lisse de dimension $n+\ell(R)-2$ et
$Z\longrightarrow Y$ est un fibr\'e en droites projectives
pour la topologie \'etale ([M1] p. 603). Le morphisme universel 
$G$-\'equivariant 
$\mathbb{P}^{1}\times V\longrightarrow X$ 
(l'action de $G$ sur $X$ \'etant triviale) 
est lisse ([K] II 3.5.4) et induit un morphisme propre et lisse
$Z\longrightarrow X$ de dimension 
relative $\ell(R)-2$.\\
\indent Montrons que $\ell(R)=2$. Soit $v\in V$. Ecrivons 
$v^{*}T_{X}=\mathcal{O}_{\mathbb{P}^{1}}(a_{1})
\oplus\cdots\oplus\mathcal{O}_{\mathbb{P}^{1}}(a_{n})$ avec 
$a_{1}\ge\cdots\ge a_{n}\ge 0$ et $a_{1}\ge 2$. Il existe un ouvert
$U\subset V$ non vide tel que le $n$-uplet 
$(a_{1},\cdots,a_{n})$ soit ind\'ependant de $v\in U$. 
Consid\'erons le morphisme :
\begin{equation*}
\begin{CD}
V @)\psi)) X\times X \\
v @))) (v(0),v(\infty))
\end{CD}
\end{equation*}
La diff\'erentielle de $\psi$ est donn\'ee par la formule ([K] II 3.4) :
\begin{equation*}
\begin{CD}
H^{0}(\mathbb{P}^{1},v^{*}T_{X}) @){d\psi(v)}))
v^{*}T_{X}\otimes k(0) \oplus v^{*}T_{X}\otimes k(\infty)\\
s @))) (s(0),s(\infty))
\end{CD}
\end{equation*}
Calculons le rang de ladite diff\'erentielle. Consid\'erons
les applications lin\'eaires $(d\psi)_{i}(v)$ ($1\le i\le n$) :
\begin{equation*}
\begin{CD}
H^{0}(\mathbb{P}^{1},\mathcal{O}_{\mathbb{P}^{1}}(a_{i})) 
@)))
\mathcal{O}_{\mathbb{P}^{1}}(a_{i})\otimes k(0) 
\oplus\mathcal{O}_{\mathbb{P}^{1}}(a_{i})\otimes k(\infty)\\
s_{i} @))) (s_{i}(0),s_{i}(\infty))
\end{CD}
\end{equation*}
Le rang de $(d\psi(v))_{i}$ est $2$ si $a_{i}\ge 1$ et $1$ si $a_{i}=0$,
de sorte que, pour $v\in U$ :
$$\text{rang}(d\psi(v))=2\text{Card}\{i|a_{i}\ge 1\}+
\text{Card}\{i|a_{i}=0\}.$$
\indent Evaluons la dimension de l'image de $\psi$. Soient $p$ et $q$
les projections de $X\times X$ sur chacun des facteurs. Le morphisme 
$Z\longrightarrow X$ \'etant propre et lisse, $p(\text{Im}(\psi))=X$. 
Si $x\in p(\text{Im}(\psi))=X$,
$p^{-1}(x)\cap\text{Im}(\psi)$ s'identifie, via la projection $q$, 
au lieu des points de $X$ par lesquels
il passe une courbe rationnelle $v(\mathbb{P}^{1})$ ($v\in V$) 
contenant $x$ et sa dimension 
est donc au moins $\ell(R)-1$ ([W] 1.11). Il en r\'esulte que la dimension 
de l'image de $\psi$ est au moins \'egale \`a $n+\ell(R)-1$. On a donc :
$$\text{dim}(\text{Im}(\psi))=2\text{Card}\{i|a_{i}\ge 1\}+
\text{Card}\{i|a_{i}=0\}\ge n+\ell(R)-1=n-1+\sum_{i=1}^{n}a_{i}.$$
Comme $a_{1}\ge 2$ :
$$n-1+\sum_{i=1}^{n}a_{i}\ge n-1+2+(\text{Card}\{i|a_{i}\ge 1\}-1)
=2\text{Card}\{i|a_{i}\ge 1\}+\text{Card}\{i|a_{i}=0\}.$$
D'o\`u :
$$\text{Card}\{i|a_{i}\ge 1\}=\ell(R)-1=\sum_{i=1}^{n}a_{i}-1$$
On en d\'eduit $\ell(R)=2$ et $(a_{1},\cdots,a_{n})=(2,0,\cdots,0)$
puisque $a_{i}=0$ ou bien $a_{i}\ge 2$
(corollaire 1). Le morphisme $Z\longrightarrow X$ est dont un
rev\^etement \'etale fini, ce qui termine la preuve du lemme.\\
\newline
\textbf{Corollaire 3}.$-$\textit{Soit $X$ une vari\'et\'e projective 
lisse de dimension $n\ge 1$ dont le fibr\'e tangent est totalement 
d\'ecompos\'e. Alors $X$ est unir\'egl\'ee si et seulement si
$X$ n'est pas minimale.}\\
\newline
\textit{D\'emonstration}.$-$ Si $X$ est unir\'egl\'ee alors $X$ n'est pas
minimale par le corollaire 2. Supposons inversement $X$ non minimale ; 
le lemme 2 entra\^{\i}nent l'existence d'un 
rev\^etement \'etale fini $Z$ de $X$, avec $Z$ 
unir\'egl\'ee. La vari\'et\'e $X$ est donc unir\'egl\'ee, ce qui termine 
la preuve du corollaire.\\
\newline
\indent Soient $X$ et $Y$ deux vari\'et\'es projectives lisses et
$X\overset{\phi}{\longrightarrow} Y$ un morphisme lisse. Une 
\textit{connexion} sur $\phi$ est un scindage de la suite exacte 
$$0\longrightarrow\phi^{*}\Omega_{Y}^{1}\longrightarrow\Omega_{X}^{1}
\longrightarrow\Omega_{X/Y}^{1}\longrightarrow 0,$$
c'est-\`a-dire un sous-fibr\'e $E\subset\Omega_{X}^{1}$ tel que 
le morphisme $E\longrightarrow\Omega_{X/Y}^{1}$ induit par la
projection $\Omega_{X}^{1}\longrightarrow\Omega_{X/Y}^{1}$ soit un
isomorphisme. 
La connexion est dite \textit{int\'egrable} si $dE\subset
E\wedge\Omega_{X}^{1}$ ou bien, ce qui est \'equivalent, si le noyau de la
surjection $T_{X}\longrightarrow E^{*}$ est int\'egrable au sens
usuel. Le morphisme $\phi$ est alors analytiquement
localement trivial de fibre $F$ et $X$
s'identifie au quotient 
de $\tilde{Y}\times F$ par le groupe $\pi_{1}(Y)$ agissant 
diagonalement o\`u $\tilde{Y}$ est le rev\^etement 
universel de $Y$ ; le scindage 
$\Omega_{X}^{1}=\phi^{*}\Omega_{Y}^{1}\oplus E$ se rel\`eve en la
d\'ecomposition 
$\Omega_{\tilde{Y}\times
F}^{1}=\Omega_{\tilde{Y}}^{1}\oplus\Omega_{F}^{1}$ ([B1] 4.5).\\
\newline
\textbf{Proposition 1}.$-$\textit{Si le th\'eor\`eme 1 est vrai en dimension 
$n\ge 1$, il est vrai en dimension $n+1$ pour les vari\'et\'es non
minimales.}\\
\newline
\textit{D\'emonstration}.$-$Soient $X$ une vari\'et\'e projective lisse non 
minimale de dimension $n+1$ ($n\ge1$)
dont le fibr\'e cotangent est totalement d\'ecompos\'e et 
$\Omega_{X}^{1}=L_{1}\oplus\cdots\oplus L_{n+1}$ ladite d\'ecomposition.
Quitte \`a passer \`a un rev\^etement \'etale, 
on peut toujours supposer
qu'il existe une vari\'et\'e projective lisse $Y$ de dimension 
$n$ et un morphisme 
$X\overset{\phi}{\longrightarrow} Y$ dont les fibres sont des droites 
projectives (lemme 2). Soit $C$ une fibre de $\phi$. 
Posons $d_{i}=L_{i}.C$ et supposons $d_{1}\ge d_{2}\ge\cdots\ge d_{n+1}$.
Consid\'erons la suite exacte :
$$0\longrightarrow N_{C/X}^{*}=\mathcal{O}_{C}^{\oplus n}
\longrightarrow{\Omega_{X}^{1}}_{|C}\longrightarrow
\Omega_{C}^{1}=\mathcal{O}_{C}(-2)\longrightarrow 0.$$
Le groupe 
$\text{Ext}^{1}_{\mathcal{O}_{C}}(\mathcal{O}_{C}(-2),
\mathcal{O}_{C}^{\oplus n})$ est nul et
l'extension pr\'ec\'edente est donc triviale. On en d\'eduit
$d_{1}=-2$ et $d_{i}=0$ pour $i\ge 2$ puis, par les th\'eor\`emes de 
changement de base, qu'il existe des fibr\'es
inversibles ${(E_{i})}_{2\le i\le n+1}$ sur $Y$ tels que 
$L_{i}=\phi^{*}E_{i}$. On v\'erifie par restriction aux fibres que
l'application $\phi^{*}\Omega_{Y}^{1}\longrightarrow L_{1}$
est identiquement nulle ; les sous-fibr\'es $\phi^{*}\Omega_{Y}^{1}$
et $L_{2}\oplus\cdots\oplus L_{n+1}$ de $\Omega_{X}^{1}$ sont donc
\'egaux et l'application 
$L_{1}\longrightarrow\Omega_{X/Y}^{1}$ est un isomorphisme. On en d\'eduit que 
$T_{Y}$ est totalement d\'ecompos\'e et on v\'erifie que la condition  
d'int\'egrabilit\'e du th\'eor\`eme est satisfaite. Le fibr\'e 
$L_{1}$ d\'etermine une connexion int\'egrable
sur $\phi$ et l'assertion souhait\'ee en r\'esulte facilement.\\
\newline
\indent La suite de notre travail est consacr\'ee \`a l'\'etude des 
vari\'et\'es minimales dont le fibr\'e tangent est totalement
d\'ecompos\'e.\\
\newline
\textbf{Lemme 3}.$-$\textit{Soient $X$ une vari\'et\'e projective lisse
minimale dont le fibr\'e tangent est totalement d\'ecompos\'e
et $T_{X}=M_{1}\oplus\cdots\oplus M_{n}$ ($n\ge 1$) ladite 
d\'ecomposition. Soit $H$ une 
section hyperplane de $X$. Alors $M_{i}H^{n-1}\le 0$ et si 
l'in\'egalit\'e pr\'ec\'edente est une \'egalit\'e on a 
$c_{1}(M_{i})=0$.}\\
\newline
\textit{D\'emonstration}.$-$La vari\'et\'e $X$ n'\'etant pas 
unir\'egl\'ee (corollaire 3) il existe une courbe lisse 
$C\in |mH|\cap\cdots\cap |mH|$ 
($m\gg 0$) telle que  
le fibr\'e ${\Omega_{X}^{1}}_{|C}$ soit num\'eriquement 
effectif ([Mi]), c'est-\`a-dire telle que le faisceau inversible  
$\mathcal{O}_{Z}(1)$ soit num\'eriquement effectif, 
o\`u $Z=\mathbb{P}_{C}({\Omega_{X}^{1}}_{|C})$.
Soit $\sigma_{i}$ la section de $Z\longrightarrow C$ associ\'ee \`a
la surjection ${\Omega_{X}^{1}}_{|C}\twoheadrightarrow {M_{i}^{-1}}_{|C}$.
On a donc 
$\sigma_{i}(C).\mathcal{O}_{Z}(1)=C.M_{i}^{-1}\ge 0$, 
ce qui prouve la 
premi\`ere assertion. La seconde r\'esulte du th\'eor\`eme de 
l'indice de Hodge puisque $M_{i}^{2}=0$ (lemme 1).\\
\newline
\textbf{Proposition 2 }([B1] th\'eor\`eme C).$-$\textit{Soient $X$ 
une vari\'et\'e projective 
lisse de dimension $n\ge1$ dont le fibr\'e tangent est 
totalement d\'ecompos\'e et  
$T_{X}=M_{1}\oplus\cdots\oplus M_{n}$ ladite d\'ecomposition.
On suppose qu'il existe une section hyperplane $H$ de $X$ telle 
que $M_{i}.H^{n-1}<0$ pour $1\le i\le n$. Le rev\^etement universel 
$\tilde{X}$ de $X$ est alors isomorphe au produit $D^{n}$ o\`u 
$D=\{z\in \mathbb{C}\quad\vert\quad |z|<1\}$ et la
d\'ecomposition de $T_{X}$ est induite par la d\'ecomposition
canonique de $T_{\tilde{X}}$.}\\
\newline
\textbf{Lemme 4}.$-$\textit{Soit $X$ une vari\'et\'e projective lisse
minimale dont le fibr\'e tangent est totalement d\'ecompos\'e. La
composante neutre $G$ du groupe des automorphismes de $X$ est alors une
vari\'et\'e ab\'elienne.}\\
\newline
\textit{D\'emonstration}.$-$Soit $T_{X}=M_{1}\oplus\cdots\oplus M_{n}$
($n\ge 1$) la d\'ecomposition de $T_{X}$. On a $H^{0}(X,M_{i})=(0)$ sauf
si le fibr\'e $M_{i}$ est trivial (lemme 3). Soit 
$k=h^{0}(X,T_{X})$ la dimension de $G$ et soit
$\mathcal{O}(x)$ une orbite de dimension minimale 
($x\in X$) en particulier ferm\'ee. Consid\'erons les applications 
tangentes :
$$T_{G,id}=H^{0}(X,T_{X})\longrightarrow T_{\mathcal{O}(x),x}
\longrightarrow T_{X,x}.$$
L'application $H^{0}(X,T_{X})\longrightarrow T_{X,x}$ obtenue est 
l'\'evaluation en $x$ qui est injective. 
On en d\'eduit que l'application
$T_{G,id}\longrightarrow T_{\mathcal{O}(x),x}$ est un isomorphisme
puis que le morphisme 
$G\longrightarrow \mathcal{O}(x)$ est un rev\^etement \'etale fini.
Les vari\'et\'es $G$ et $\mathcal{O}(x)$ 
sont donc des vari\'et\'es ab\'eliennes de dimension $k$.\\
\newline
\textbf{Lemme 5 }([B2] proposition 2.1).$-$\textit{Soient $X$ 
une vari\'et\'e projective lisse de dimension $n\ge2$ et 
$M$ un fibr\'e inversible non trivial 
dont le premi\`ere classe de Chern $c_{1}(M)\in H^{2}(X,\mathbb{Z})$ 
est nulle. On suppose qu'il existe deux 1-formes non nulles 
$\alpha\in H^{0}(X,\Omega_{X}^{1}\otimes M)$ et 
$\omega\in H^{0}(X,\Omega_{X}^{1})$ telles que 
$\alpha\wedge\omega=0$. Alors il existe un morphisme 
surjectif \`a fibres connexes 
$X\longrightarrow C$ vers une courbe lisse $C$ de genre 
$g(C)\ge1$ tel que $\omega$ provienne par image r\'eciproque de 
$H^{0}(C,\Omega_{C}^{1})$.}\\
\newline
\indent Soient $X$ une vari\'et\'e projective lisse et 
$X\overset{\phi}{\longrightarrow} C$ un morphisme surjectif vers une courbe 
lisse $C$. Le \textit{diviseur de ramification} de $\phi$ est d\'efini par 
la formule :
$$D(\phi)=\sum_{p\in C}\phi^{*}p-(\phi^{*}p)_{red}\,;$$
on d\'emontre que les sections du fibr\'e 
$\phi^{*}\omega_{C}(D(\phi))$ sont holomorphes, \textit{i.e.}, 
$\phi^{*}\omega_{C}(D(\phi))\subset\Omega_{X}^{1}$ ([D] lemme 4.4).\\ 
\newline
\textbf{Lemme 6 }([B2]).$-$\textit{Soient $X$ une vari\'et\'e projective lisse 
de dimension $n\ge 2$, 
$X\overset{\phi}{\longrightarrow} C$ un morphisme surjectif vers une courbe 
projective lisse $C$ et $D$ un diviseur vertical. On a $D^{2}\le 0$
et si $D^{2}=0$ alors il existe $r\in\mathbb{Q}^{*}$ tel 
que $rD\in\phi^{*}\text{Pic}(C)$.}\\
\newline
\textbf{Lemme 7}.$-$\textit{Soient $X$ une vari\'et\'e projective lisse 
de dimension $n\ge 2$ et 
$X\overset{\phi}{\longrightarrow} C$ un morphisme surjectif vers une courbe 
lisse $C$. Soit 
$L=\phi^{*}\omega_{C}(D(\phi))\subset\Omega^{1}_{X}$ o\`u $D(\phi)$ 
est le diviseur de ramification de $\phi$. Si $L$
est un sous-fibr\'e du fibr\'e cotangent et $L^{2}=0$ alors
les fibres de $\phi$ sont lisses ou 
multiples de vari\'et\'es lisses.}\\
\newline
\textit{D\'emonstration}.$-$ Ecrivons 
$D(\phi)=D_{1}+\cdots+D_{k}$ ($k\ge1$) où
$D_{i}=\phi^{*}p_{i}-{(\phi^{*}p_{i})}_{red}$, les points $p_{i}\in C$
étant tous distincts et soit $\phi^{*}p_{i}=\sum_{j}n_{i,j}D_{i,j}$ la
décomposition de la fibre schématique $\phi^{*}p_{i}$ en somme de ses
composantes irréductibles ($n_{i,j}\ge1$). Le diviseur $D(\phi)$ 
est vertical et on a donc ${D_{i}}^{2}=0$ pour tout entier 
$i\in\{1,\cdots,k\}$ (lemme 6). On en d\'eduit qu'il existe 
$r_{i}\in\mathbb{Q}^{*}$ tel 
que $r_{i}D_{i}\in\phi^{*}\text{Pic}(C)$ (lemme 6).
On a donc $n_{i,j}=n_{i,k}=n_{i}$ pour tout triplet
$(i,j,k)$ et
$D_{i}=(n_{i}-1){(\phi^{*}p_{i})}_{red}$. Prenons un point $p_{i}\in C$ 
tel que la fibre schématique
$\phi^{*}p_{i}$ soit non réduite. Soient $V\subset C$ un disque
ouvert pour la topologie usuelle contenant le point $p_{i}$ et 
$U=\phi^{-1}(V)$. On suppose que $p_{i}$ est le seul point de
$V$ au dessus duquel la fibre est non réduite. En effectuant
le changement de base local $z\mapsto z^{n_{i}}$, puis en normalisant, 
on obtient un diagramme commutatif :  
\begin{equation*}
\begin{CD}
U_{1} @)q))U \\
@VV\phi_{1}V @VV{\phi}V\\
V_{1} @)z\mapsto z^{n_{i}})) V
\end{CD}
\end{equation*}
\noindent où $q$ est un revêtement étale fini et
$\phi_{1}$ est un morphisme à fibres réduites et connexes. 
L'inclusion  
$L\subset\Omega_{X}^{1}$ se rel\`eve en 
l'inclusion naturelle
$\phi_{1}^{*}\omega_{V_{1}}\subset
\Omega_{U_{1}}^{1}$ et le morphisme $\phi_{1}$ est donc 
lisse, ce qui termine la preuve du lemme.\\
\newline
\textbf{Proposition 3}.$-$\textit{Soient $X$ une vari\'et\'e
projective lisse minimale de dimension $n\ge1$ dont le fibr\'e tangent
est totalement d\'ecompos\'e et $M$ un facteur direct de ladite
d\'ecomposition. Si $c_{1}(M)=0$ alors $M$ est de torsion.}\\
\newline
\textit{D\'emonstration}.$-$ Soient $\Omega_{X}^{1}=L_{1}\oplus\cdots\oplus
L_{n}$ la d\'ecomposition du fibr\'e cotangent et $H$ une section
hyperplane de $X$. On suppose
que l'un des facteurs $L_{i}$ $(1\le i\le n)$ 
est de premi\`ere classe de Chern 
nulle mais pas de torsion.\\
\indent Un fibr\'e inversible de torsion est trivialis\'e
par un rev\^etement \'etale fini et un fibr\'e inversible qui n'est
pas de torsion le reste apr\`es ledit rev\^etement. On peut donc
supposer que les
fibr\'es $L_{1},\cdots,L_{k}$ sont de premi\`ere classe de Chern nulles mais pas
de torsion ($k\ge 1$) et que les fibr\'es $L_{k+1},\cdots,L_{k+r}$
sont triviaux. 
On a enfin $L_{i}H^{n-1}>0$ pour $i\ge k+r+1$ (lemme 3). On a donc 
$h^{0}(X,T_{X})=r$ et pour tout rev\^etement \'etale fini
$Z\longrightarrow X$ on a $h^{0}(Z,T_{Z})=r$.\\
\indent La th\'eorie de Hodge fournit un
isomorphisme antilin\'eaire 
$H^{1}(X,L_{1})\simeq H^{0}(X,\Omega^{1}_{X}\otimes L_{1}^{-1})$ 
([B2] 3.4) puisque $c_{1}(L_{1})=0$ et on a donc
$H^{1}(X,L_{1})\neq0$. Il existe alors 
deux 1-formes non nulles 
$\alpha\in H^{0}(X,\Omega_{X}^{1}\otimes L_{1}^{-1})$
et $\omega\in H^{0}(X,\Omega_{X}^{1})$ telles que 
$\alpha\wedge\omega=0$ ([B2] th\'eor\`eme 2.2 et [S]) et un 
fibr\'e inversible $L\subset\Omega_{X}^{1}$ tel que 
$\omega\in H^{0}(X,L)$ (lemme 5). La 1-forme $\alpha$ d\'etermine 
une application injective $L_{1}\longrightarrow\Omega_{X}^{1}$ dont  
l'image sera not\'ee $M$. La condition $\alpha\wedge\omega=0$ signifie 
que les deux sous-faisceaux $L$ et $M$ de $\Omega_{X}^{1}$ coincident 
sur un ouvert non vide de $X$. Or il existe un entier 
$m\in\{1,\cdots,n\}$ tel que la projection 
$M\longrightarrow L_{m}$ soit non nulle. La composante 
$\alpha_{m}$ de $\alpha$ sur le facteur direct 
$H^{0}(X,L_{m}\otimes L_{1}^{-1})$ de 
$H^{0}(X,\Omega_{X}^{1}\otimes L_{1}^{-1})$ est 
donc non nulle. La projection $L\longrightarrow L_{m}$ est 
\'egalement non nulle puisque les sous-faisceaux $L$ et $M$ de 
$\Omega_{X}^{1}$ coincident g\'en\'eriquement. Notons 
$\omega_{m}$ la composante non nulle de $\omega$ sur le facteur direct 
$H^{0}(X,L_{m})$ de $H^{0}(X,\Omega_{X}^{1})$. Par choix des 
formes $\alpha_{m}$ et $\omega_{m}$ on a 
$\alpha_{m}\wedge\omega_{m}=0$. On en d\'eduit qu'il existe 
un morphisme surjectif 
\`a fibres connexes $X\overset{\phi}{\longrightarrow} C$ vers une courbe 
$C$ de genre $g(C)\ge 1$ tel que la forme 
$\omega_{m}$ provienne par image r\'eciproque de 
$H^{0}(C,\Omega_{C}^{1})$ (lemme 5). Les sous-faisceaux
$\phi^{*}\omega_{C}$ et 
$L_{m}$ de $\Omega_{X}^{1}$ coincident sur un ouvert non vide
puisque $\omega_{m}\in H^{0}(X,L_{m})$. 
On en d\'eduit 
l'\'egalit\'e $L_{m}=\phi^{*}\omega_{C}(D(\phi))$ o\`u $D(\phi)$ 
est le
diviseur de ramification de $\phi$ ([D] lemmes 4.1, 4.2 et 4.4).
Or $L_{m}^{2}=0$ (lemme 1) et 
les fibres de $\phi$ sont donc lisses ou multiples de 
vari\'et\'es lisses (lemme 7).
La courbe $C$ \'etant de genre au moins $1$, il existe 
un morphisme fini $C_{1}\overset{\pi}{\longrightarrow} C$ tel
que l'indice de ramification de $\pi$ en $q\in C_{1}$ 
soit \'egal \`a la multiplicit\'e de la fibre $\phi^{-1}(\pi(q))$ 
([K-O] lemme 6.1). 
Soit $X_{1}$ la normalisation du produit fibr\'e
$X{\times}_{C} C_{1}$ et soient 
$X_{1}\overset{\phi_{1}}{\longrightarrow}C_{1}$ 
et
$X_{1}\overset{\pi_{1}}{\longrightarrow} X$
les morphismes natuels. Le morphisme $\phi_{1}$ 
est lisse et $\pi_{1}$ est un 
rev\^etement \'etale. Les sous-fibr\'es  
$\phi_{1}^{*}\omega_{C_{1}}$ et $\pi_{1}^{*}L_{m}$
de $\Omega_{X_{1}}^{1}$ coincident et le fibr\'e 
${\oplus}_{i\neq m}\pi_{1}^{*}L_{i}$ d\'etermine donc une connexion 
sur $\phi_{1}$ qui est automatiquement int\'egrable puisque $C_{1}$
est une courbe. Le morphisme $\phi_{1}$ est donc analytiquement localement trivial
de fibre $F$. De plus il existe un morphisme non nul
$\pi_{1}^{*}L_{1}\longrightarrow \pi_{1}^{*}L_{m}$. On a donc
$\pi_{1}^{*}L_{m}=\pi_{1}^{*}L_{1}(D)$ o\`u $D$ est un
diviseur effectif. On en d\'eduit que $D$ est vertical puisque
$c_{1}(\pi_{1}^{*}L_{1})=0$ et 
$\pi_{1}^{*}L_{m}=\phi_{1}^{*}\omega_{C_{1}}$ puis qu'il existe
un fibr\'e $M_{1}$ sur $C_{1}$ tel que
$\pi_{1}^{*}L_{1}=\phi_{1}^{*}M_{1}$. Remarquons enfin que
$g(C_{1})\ge 2$. Sinon $C_{1}$ serait une courbe elliptique et
puisqu'on a une application non triviale
$\pi_{1}^{*}L_{1}\longrightarrow
\pi_{1}^{*}L_{m}=\phi_{1}^{*}\omega_{C_{1}}$ le fibr\'e $L_{1}^{-1}$ serait
donc trivial ce qui est contraire aux hypoth\`eses. On 
a donc $m\neq 1$ et $c_{1}(\pi_{1}^{*}L_{m})\neq 0$.\\ 
\indent Le fibr\'e tangent $T_{F}$ est totalement
d\'ecompos\'e et $F$ est minimale puisque $K_{F}={K_{X_{1}}}_{|F}$. La
composante neutre $G$ du groupe des automorphismes de $F$ est donc une
vari\'et\'e ab\'elienne (lemme 4).
On a enfin $h^{0}(F,T_{F})\ge h^{0}(X_{1},T_{X_{1}})+1$
puisque la restriction du facteur $\pi_{1}^{*}L_{1}$ \`a $F$ est
triviale. Consid\'erons le
sch\'ema en groupes $\text{Aut}^{0}(X_{1}/C_{1})$ au dessus de $C_{1}$
dont la fibre au dessus d'un point $p\in C_{1}$ est la composante
neutre du groupe des automorphismes de $\phi_{1}^{-1}(p)$.
Ladite fibre est une vari\'et\'e ab\'elienne
isomorphe \`a $G$ puisque 
les fibres de $\phi_{1}$ sont toutes 
isomorphes \`a $F$. Il existe donc un rev\^etement \'etale fini 
$C_{2}\longrightarrow C_{1}$ trivialisant ledit sch\'ema en
groupes. Cette assertion r\'esulte de l'existence d'un espace de
modules fin pour les vari\'et\'es ab\'eliennes munies d'une
polarisation de degr\'e donn\'e et d'une structure de niveau
$n\ge3$. On en d\'eduit que le groupe $G$ agit sur le produit fibr\'e  
$X_{2}=X_{1}\times_{C_{1}}C_{2}$ puis que
$\text{dim}\,G=h^{0}(F,T_{F})\le
h^{0}(X_{2},T_{X_{2}})$. Or $h^{0}(X_{2},T_{X_{2}})=
h^{1}(X_{1},T_{X_{1}})\le h^{0}(F,T_{F})-1$ 
puisque la projection $X_{2}\longrightarrow X_{1}$ est un rev\^etement
\'etale, ce qui donne la contradiction cherch\'ee.\\
\newline
\textbf{Lemme 8}.$-$\textit{Soit $X$ une vari\'et\'e projective lisse
de dimension $n\ge 1$. On suppose que la composante neutre $G$ du groupe des 
automorphismes de $X$ est une vari\'et\'e ab\'elienne de dimension $k\ge 1$
et que l'application canonique 
$H^{0}(X,T_{X})\longrightarrow H^{0}(X,\Omega_{X}^{1})^{*}$
est injective. Il existe alors un rev\^etement \'etale fini 
$G$-\'equivariant $G\times F\longrightarrow X$, $G$ agissant 
diagonalement sur $G\times F$ par translations sur $G$ et trivialement
sur $F$.}\\
\newline
\textit{D\'emonstration}.$-$ Soient  
$A=H^{0}(X,\Omega_{X}^{1})^{*}/\text{Im}(H_{1}(X,\mathbb{Z}))$ 
la vari\'et\'e d'Albanese de $X$ et $x\in X$ un point de $X$. 
Consid\'erons le morphisme d'Albanese :
\begin{equation*}
\begin{CD}
X @)a)) A \\
y @))) [\omega\mapsto\int_{x}^{y}\omega]
\end{CD}
\end{equation*}
Le groupe $G$ \'etant connexe,  
la repr\'esentation canonique 
$G\longrightarrow GL(H^{0}(X,\Omega_{X}^{1}))$ est triviale et 
toute 1-forme holomorphe sur $X$ est donc $G$-invariante. On en d\'eduit
que pour tout couple $(y,z)\in X\times X$ et tout $g\in G$ on a  
$[\omega\mapsto\int_{y}^{g(y)}\omega]
= 
[\omega\mapsto\int_{z}^{g(z)}\omega]$ dans $A$. Notons 
$\mathcal{O}(x)$ l'orbite de $x$ sous $G$.
L'application :
\begin{equation*}
\begin{CD}
G @))) A \\
g @))) [\omega\mapsto\int_{x}^{g(x)}\omega]
\end{CD}
\end{equation*}
obtenue par composition de l'inclusion naturelle 
$\mathcal{O}(x)\subset X$ et du morphisme d'Albanese 
est donc un morphisme de groupes alg\'ebriques ind\'ependant 
du point $x\in X$ consid\'er\'e et le morphisme d'Albanese 
est $G$-\'equivariant, $G$ agissant par translations sur $A$ 
via le morphisme $G\longrightarrow A$.
Notons que l'image 
de $G\longrightarrow A$ est une vari\'et\'e ab\'elienne de dimension $k$.
En effet, l'application tangente dudit morphisme est 
l'application canonique
$H^{0}(X,T_{X})\otimes\mathcal{O}_{G}\longrightarrow 
H^{0}(X,\Omega_{X}^{1})^{*}\otimes\mathcal{O}_{G}$, qui est injective par
hypoth\`ese. On en 
d\'eduit qu'il existe une sous-vari\'et\'e ab\'elienne $B\subset A$ de 
dimension $\text{dim}(A)-k$ telle que le morphisme canonique 
$G\times B\longrightarrow A$ soit \'etale fini 
([Mu] th\'eor\`eme de Poincar\'e).
Soient $X_{1}$ le produit fibr\'e $X\times_{A}(G\times B)$ 
et $q$ la projection sur $G$. Le groupe $G$ agit diagonalement 
sur le produit $G\times B$ en agissant par translations sur lui-m\^eme 
et trivialement sur $B$ ; le morphisme 
$G\times B\longrightarrow A$ est alors $G$-\'equivariant. On en d\'eduit  
que $G$ agit librement sur $X_{1}$ et que le morphisme  
$q$ est $G$-\'equivariant.
Il en r\'esulte que $X_{1}$ est isomorphe, au dessus de $G$,
au produit $G\times F$, $F$ \'etant une fibre de $q$, ce qui termine
la preuve du lemme.\\
\newline
\textbf{Lemme 9}.$-$\textit{Soient $G$ une vari\'et\'e ab\'elienne de
dimension $n\ge1$ et
$F$ une vari\'et\'e projective lisse de dimension $m\ge1$. Alors toute
connexion sur la projection $G\times F\longrightarrow F$ est
int\'egrable.}\\
\newline
\textit{D\'emonstration}.$-$Notons $E\subset\Omega_{G\times F}^{1}$ ladite
connexion. Le fibr\'e $E$ est trivial et il existe donc des $1$-formes
$\omega_{i}\in H^{0}(G\times F,\Omega_{G\times F}^{1})$ $(1\le i\le
n)$ telles que $E=\mathcal{O}_{G\times
F}\omega_{1}\oplus\cdots\oplus\mathcal{O}_{G\times F}\omega_{n}$. Soient
$f_{1},\cdots,f_{n}$ des fonctions holomorphes d\'efinies sur un
ouvert de $G\times F$. On a la formule :
$$d(\sum_{i=1}^{n}f_{i}\omega_{i})=\sum_{i=1}^{n}df_{i}\wedge\omega_{i}+
\sum_{i=1}^{n}f_{i}d\omega_{i}.$$
Or $d\omega_{i}=0$ pour $i\in\{1,\cdots,n\}$ puisque $G\times F$ est
projective. On a donc
$dE\subset E\wedge\Omega_{G\times F}^{1}$ ce qui termine la preuve du
lemme.\\
\newline
\textbf{Th\'eor\`eme 2}.$-$\textit{Soit X une vari\'et\'e projective
lisse minimale de dimension $n\ge 1$ dont le fibr\'e tangent est totalement 
d\'ecompos\'e. Le rev\^etement universel 
$\tilde{X}$ de $X$ est alors produit de surfaces de Riemann 
et la d\'ecomposition de $T_{X}$ est induite 
par la d\'ecomposition canonique de $T_{\tilde{X}}$.}\\
\newline
\textit{D\'emonstration}.$-$Soit
$\Omega_{X}^{1}=L_{1}\oplus\cdots\oplus L_{n}$ $(n\ge 1)$ la 
d\'ecomposition du fibr\'e cotangent et soit $H$ une section hyperplane de $X$.
On a $c_{1}(L_{i})H^{n-1}\ge0$ (lemme 3) et si $c_{1}(L_{i})H^{n-1}>0$ pour tout 
$i\in\{1,\cdots,n\}$, la proposition 2 permet de conclure. Quitte \`a
permuter les $L_{i}$ on peut toujours supposer que 
les fibr\'es $L_{1},\cdots,L_{k}$ ($k\ge1$) sont de torsion et que 
$c_{1}(L_{i})\neq 0$ pour $i\ge k+1$ (proposition 3). Il existe un
rev\^etement \'etale fini de $X$ trivialisant ces fibr\'es et 
il suffit donc de traiter le cas o\`u lesdits fibr\'es 
sont triviaux. On a alors $h^{0}(X,L_{i}^{-1})=0$ pour $i\ge k+1$
(lemme 3) et $h^{0}(X,T_{X})=k$. Soit $G$ la composante
neutre du groupe des automorphismes de $X$ ; $G$ est une vari\'et\'e
ab\'elienne de dimension $k$ (lemme 4). L'application canonique
$H^{0}(X,T_{X})\longrightarrow 
H^{0}(X,\Omega_{X}^{1})^{*}$ est injective et il existe donc un
rev\^etement \'etale fini 
$G\times F\overset{p}{\longrightarrow} X$ (lemme 8) o\`u $F$ est 
une vari\'et\'e projective lisse. Le morphisme
$p$ \'etant \'etale,
on a 
$\Omega_{G\times F}^{1}=p^{*}\Omega_{X}^{1}=p^{*}L_{1}\oplus\cdots\oplus p^{*}L_{n}$ 
et 
$(p^{*}L_{i})({p^{*}H}^{n-1})>0$ pour $i\ge k+1$ (lemme 3). Par hypoth\`ese, les 
fibr\'es $p^{*}L_{i}$ sont triviaux pour $i\in\{1,\cdots,k\}$.
On en d\'eduit 
$h^{0}(G\times F,p^{*}L_{i}^{-1})=0$ pour $i\ge k+1$ puis
$h^{0}(G\times F,T_{G\times F})=k$ et
$h^{0}(G\times F,s^{*}T_{F})=0$ o\`u $s$ d\'esigne 
la projection
de $G\times F$ sur $F$. Notons $r$ la projection sur $G$.
L'application 
$s^{*}\Omega_{F}^{1}\longrightarrow p^{*}L_{1}\oplus\cdots\oplus p^{*}L_{k}$ 
est identiquement nulle et les sous-fibr\'es $s^{*}\Omega_{F}^{1}$ et 
$p^{*}L_{k+1}\oplus\cdots\oplus p^{*}L_{n}$ 
de $\Omega_{G\times F}^{1}$ sont donc isomorphes. On en d\'eduit que  
l'application  
$p^{*}L_{1}\oplus\cdots\oplus p^{*}L_{k}\longrightarrow r^{*}T_{G}$ est 
un isomorphisme. Soit $x$ un point de $F$. Le fibr\'e
${p^{*}L_{i}}_{|G\times\{x\}}$ $(k+1\le i\le n)$ est donc facteur 
direct du fibr\'e trivial 
$\mathcal{O}_{G\times\{x\}}\oplus\cdots\oplus\mathcal{O}_{G\times\{x\}}$
et donc 
lui-m\^eme trivial. Il existe donc des fibr\'es
$E_{i}\in\text{Pic}(F)$ $(k+1\le 
i\le n)$ tels que $p^{*}L_{i}=s^{*}E_{i}$. Soit $H_{F}$ une
section hyperplane de $F$. On a donc $E_{i}H_{F}^{n-k}>0$ (lemme 3)
puisque $c_{1}(L_{i})\neq 0$ pour $i\ge k+1$. On a
$\Omega_{F}^{1}=E_{k+1}\oplus\cdots\oplus E_{n+1}$ et la proposition 2
permet de conclure puisque le fibr\'e 
$p^{*}L_{1}\oplus\cdots\oplus p^{*}L_{k}$ d\'etermine une
connexion int\'egrable sur la projection $s$ (lemme 9).\\
\newline  
\textbf{Remarque }([B2]).$-$L'hypoth\`ese d'int\'egrabilit\'e dans le
th\'eor\`eme 1 
est effectivement n\'ec\'essaire. En effet, notons $X$ le produit
$A\times\mathbb{P}^{1}$ o\`u $A$ 
est une surface ab\'elienne. Il n'est pas difficile 
d'exhiber une connexion non int\'egrable sur la projection 
canonique $X\longrightarrow A$ ; la  
d\'ecomposition de $T_{X}$ ainsi obtenue ne se rel\`eve donc pas en la 
d\'ecomposition canonique du fibr\'e tangent $T_{\tilde{X}}$ 
o\`u $\tilde{X}$ est le rev\^etement universel de $X$.
\newpage
\centerline{\textsc{R\'ef\'erences}}
$\ $
\vspace{0.5cm}\\
\noindent [B1] A.Beauville, \emph{Complex manifolds with split tangent 
bundle}, Vol. en m\'emoire de M.Schneider, de Gruyter, \`a para\^{\i}tre.\\
\newline
\noindent [B2] A.Beauville, \emph{Annulation du $H^{1}$ pour les
fibr\'es en droites plats}, Lecture Notes in Math., 1507, 1-15, 1992.\\
\newline
\noindent [D] S.Druel, \emph{Structures de Poisson sur les vari\'et\'es
alg\'ebriques de dimension trois}, \`a para\^{\i}tre dans le Bulletin de 
la Soci\'et\'e Math. de France.\\
\newline
[K-O] S.Kobayashi, T.Ochiai, \emph{Holomorphic structures modeled
after hyperquadrics}, T\^ohoku Math. J., 34, 587-629, 1982.\\ 
\newline
\noindent [K] J.Koll\'ar, \emph{Rational curves on algebraic varieties}
Springer, Berlin$-$Heidelberg$-$New-York$-$Tokyo, 32, 1996.\\
\newline
\noindent [Mi] Y.Miyaoka, \emph{The Chern classes and Kodaira
dimension of minimal variety}, Advanced Study in Pure Math., 10,
449-476, 1987.\\
\newline
\noindent [M1] S.Mori, \emph{Projective manifolds with ample tangent 
bundles}, Ann. of Math., 110, 593-606, 1979.\\
\newline
\noindent [M2] S.Mori, \emph{Threefolds whose canonical bundles are not
numerically effective}, Ann. of Math., 116, 133-176, 1982.\\
\newline
\noindent [Mu] D.Mumford, \emph{Abelian Varieties}, Oxford University 
Press, 1970.\\
\newline
\noindent [S] C.Simpson, \emph{Subspaces of moduli spaces of rank one
local systems}, Ann. Sci. \'Ecole Norm. Sup., 26, 361-401, 1993.\\
\newline
\noindent [W] J.Wisniewski, \emph{Length of extremal rays and 
generalized adjunction}, Math. Z., 200, 409-427, 1989.\\

\end{document}